\documentclass[12pt]{article}
    \usepackage[cp866]{inputenc}
    \usepackage[russian]{babel}
    \usepackage{multicol}
    \usepackage{euscript,graphicx}
    \usepackage{amsmath,amsfonts,amssymb}
    \usepackage{euscript}
    \date{}
    \title{ }
    \author{ }
    \date{ }
    
\begin{document}
    \begin{flushright}
    УДК 514.1
    \end{flushright}
    \begin{center}
    {\bf G.G. Mihailichenko}
    \end{center}
    \begin{center}
    {\Large \bf The Phenomenologically Symmetric Geometry\\
    of Two Sets of Rank (3,2)}\\
    \end{center}

    \vspace{10mm}

    {\it Annotation}. This note deals with one of
    the simplest phenomenologically
    symmetric geometries of two sets of rank (3,2), determined on a one-dimensional
    and two-dimensional manifolds by the metric function $f = f(x,\xi,\eta)$.

    \vspace{5mm}


    \vspace{5mm}

    Suppose we are given a one-dimensional and a two-dimensional manifolds $\mathfrak{M}$ and
    $\mathfrak{N}$, whose points are defined by lowercase Latin and Greek
    letters correspondingly, and a metric
    (two-point) function $f$ that assigns to some pairs $<\!i\alpha\!> \
     \in \mathfrak{M \times N}$ a real number $f(i\alpha) \in R$. We shall assume, for one thing, that two
    axioms are satisfied as follows:

    \vspace{3mm}

    {\bf A1.} The domain of the function $f$ is open and dense in \
    $\mathfrak{M \times N}$.

    {\bf A2.} The metric function $f$ is sufficiently smooth.

    \vspace{3mm}

    If $x$ and  $\xi, \eta$ are local coordinates in the manifolds
    $\mathfrak{M}$ and $\mathfrak{N}$, then it is possible to write down for the metric function $f$
    its coordinate representation:
    $$
    f=f(x,\xi,\eta),
    \eqno(1)
    $$
    and, for example, $f(i\alpha) = f(x_i,\xi_\alpha,\eta_\alpha)$.

    Suppose one more and third axiom is satisfied:

    \vspace{3mm}

    {\bf A3.} The local coordinates \ $x$ \ and  \ $\xi, \eta$ \ are included in the
    representation (1) in an essential way.

    \vspace{3mm}

    The mathematical expression of the axiom ${\bf A.3}$ will be the following inequalities:
    $$
    \partial f(x,\xi,\eta)/\partial x\neq0, \
    \partial(f(x_1,\xi,\eta),f(x_2,\xi,\eta))/\partial(\xi,\eta)\neq0,
    \eqno(2)
    $$
    if $x_1\neq x_2$, which we shall also write for some points of the manifolds:
    $$
    \partial f(x_i,\xi_\alpha,\eta_\alpha)/\partial x_i \neq0, \
    \partial(f(x_j,\xi_\alpha,\eta_\alpha),f(x_k,\xi_\alpha,\eta\alpha))
    /\partial(\xi_\alpha,\eta_\alpha)\neq0,
    \eqno(2')
    $$
    if $j \neq k$.

    We shall construct function $F$, by assigning to the cortege
    $<ijk,\alpha\beta> \in
    \mathfrak{M^3 \times N^2}$ a point $(f(i\alpha), \\ f(i\beta),f(j\alpha),f(j\beta),
    f(k\alpha),f(k\beta)) \in R^6$, and introduce one last axiom:

    \vspace{3mm}

    {\bf A4.} The set of values of the function $F$ lies in $R^6$ on a smooth non-degenerate hyper-surface:
    $$
    \Phi(f(i\alpha),f(i\beta),f(j\alpha),f(j\beta),f(k\alpha),f(k\beta))=0.
    \eqno(3)
    $$

    {\bf Definition.} We shall say that a metric function $f$
    gives on a one-dimensional manifold $\mathfrak{M}$ and a two-dimensional manifold $\mathfrak{N}$
    {\it a phenomenologically symmetric geometry of two sets} (PS G2S)
    of rank (3,2) if all the above said axioms {\bf A1--A4} are satisfied.

    \vspace{3mm}

    The equation (3) establishes a nontrivial connection of the six functions
    $ f(i\alpha),f(i\beta),\\ f(j\alpha),f(j\beta),f(k\alpha),f(k\beta)$ of seven
    variables $x_i,x_j,x_k,\xi_\alpha,\eta_\alpha,\xi_\beta,\eta_\beta$, for which to exist
    it is necessary and sufficient that the rank of the corresponding
    functional matrix
    $$
    \left\|
    \begin{array}{cccccc}
    f_x(i\alpha ) & f_x(i\beta ) & 0 & 0 & 0 & 0 \\
    0 & 0 & f_x(j\alpha ) & f_x(j\beta ) & 0 & 0 \\
    0 & 0 & 0 & 0 & f_x(k\alpha ) & f_x(k\beta ) \\
    f_\xi(i\alpha ) & 0 & f_\xi(j\alpha) & 0 & f_\xi(k\alpha) & 0 \\
    f_\eta(i\alpha ) & 0 & f_\eta(j\alpha) & 0 & f_\eta(k\alpha) & 0 \\
    0 & f_\xi(i\beta ) & 0 & f_\xi(j\beta) & 0 & f_\xi(k\beta) \\
    0 & f_\eta(i\beta ) & 0 & f_\eta(j\beta) & 0 & f_\eta(k\beta)
    \end{array}
    \right\|
    \eqno(4)
    $$
    of dimension $7\times 6$ should be less than 6.

    \vspace{3mm}

    {\bf Theorem.} {\it If a metric function $f$
    gives on a one-dimensional manifold $\mathfrak{M}$ and a two-dimensional manifold $\mathfrak{N}$
    a phenomenologically symmetric geometry of two sets (PS G2S) of rank (3,2),
    then, with an accuracy up to a change of coordinates in the manifolds and its scaling
    transformation, the representation (1) may be written in
    the following canonical form;}
    $$
    f=x\xi+\eta.
    \eqno(5)
    $$

    \vspace{3mm}

    From the theorem, it is clear what a most general expression for the coordinate
    representation (1) of the metric function giving PS G2S of rank (3.2) must be:
    $$
    f=\chi(\varphi(x)\psi_1(\xi,\eta)+\psi_2(\xi,\eta)),
    \eqno(5')
    $$
    where $\chi,\varphi$ are non-constant functions of one variable,
    and $\psi_1,\psi_2$ are independent functions of two variables.

    Let us take two Jacobians of sixth order vanishing to zero that are obtained from the matrix (4)
    by way of crossing out the sixth and seventh rows and expand them with respect to the elements of the
    first column:
    $$
    \left.\begin{array}{c}
     - f_x(i\alpha)f_\xi(i\beta)f_x(j\beta)f_x(k\beta)A(jk,\alpha) -
    f_\xi(i\alpha)f_x(i\beta)B_1(jk,\alpha\beta) \ + \\ + \ f_\eta(i\alpha)f_x(i\beta)
    C_1(jk,\alpha\beta) = 0,   \\

     - f_x(i\alpha)f_\eta(i\beta)f_x(j\beta)f_x(k\beta)A(jk,\alpha) -
    f_\xi(i\alpha)f_x(i\beta)B_2(jk,\alpha\beta) \ + \\ + \ f_\eta(i\alpha)f_x(i\beta)
    C_2(jk,\alpha\beta) = 0,
    \end{array}\right\}
    \eqno(6)
    $$
    where
    $$
    A(jk,\alpha\beta) = \left| \begin{array}{cccc}
    f_\xi(j\alpha ) &  f_\xi(k\alpha) \\
    f_\eta(j\alpha ) &  f_\eta(k\alpha)
    \end{array} \right| \neq0,
    $$

    $$
    B_1(jk,\alpha\beta) = \left| \begin{array}{cccc}
    f_x(j\alpha ) & f_x(j\beta ) & 0 & 0 \\
    0 & 0 & f_x(k\alpha ) & f_x(k\beta ) \\
    f_\eta(j\alpha ) & 0 & f_\eta(k\alpha) & 0 \\
    0 & f_\xi(j\beta ) & 0 & f_\xi(k\beta)
    \end{array} \right|,
    $$

    $$
    C_1(jk,\alpha\beta) = \left| \begin{array}{cccc}
    f_x(j\alpha ) & f_x(j\beta ) & 0 & 0 \\
    0 & 0 & f_x(k\alpha ) & f_x(k\beta ) \\
    f_\xi(j\alpha ) & 0 & f_\xi(k\alpha) & 0 \\
    0 & f_\xi(j\beta ) & 0 & f_\xi(k\beta)
    \end{array} \right|,
    $$
    and the minors $B_2(jk,\alpha\beta)$ and $C_2(jk,\alpha\beta)$ are obtained from the minors
    $B_1(jk,\alpha\beta)$ and $C_1(jk,\alpha\beta)$ respectively by way of the change $f_\xi(j\beta) \to f_\eta(j\beta) ,f_\xi(k\beta) \to f_\eta(k\beta)$.

    As result, we get a system of two differential functional relations
    (5), the rank of the matrix of the coefficients of which, with the derivatives $f_x(i\alpha),f_\xi(i\alpha), \\ f_\eta(i\alpha)$,  is equal to two, as, for example,
    $$
    \left| \begin{array}{cccc}
    B_1(jk,\alpha\beta ) &  C_1(jk,\alpha\beta) \\
    B_2(jk,\alpha\beta) &  C_2(jk,\alpha\beta)
    \end{array} \right| \neq0.
    $$

    Indeed, expanding that determinant we do get the product
    $$
    f_x(j\alpha)f_x(j\beta)f_x(k\alpha)f_x(k\beta)
     \left| \begin{array}{cccc}
    f_\xi(j\alpha ) &  f_\xi(k\alpha) \\
    f_\eta(j\alpha ) &  f_\eta(k\alpha)
    \end{array} \right|
    \times \left| \begin{array}{cccc}
    f_\xi(j\beta ) &  f_\xi(k\beta) \\
    f_\eta(j\beta ) &  f_\eta(k\beta)
    \end{array} \right| \neq 0,
    $$
          in which each factor is unequal to zero, according to the inequalities (2).

   We shall divide each of the relations (6) by $f_x(i\beta)A(jk,\alpha)$, and then
    fix the points $j, k$ and $\beta$. Then, with respect to the function $f(i\alpha)$, we get a system of two linear homogeneous differential
    equations in partial derivatives of the first order, the rank of the system being equal to two.
    We shall write it down, introducing convenient designation for the coefficients and omitting the point indexes $i$ and $\alpha$
               of the variables (coordinates):
    $$
    \left.\begin{array}{c}
    \lambda_1f_x + \sigma_1f_\xi + \tau_1f_\eta = 0, \\
    \lambda_2f_x + \sigma_2f_\xi + \tau_2f_\eta = 0,
    \end{array}\right\}
    \eqno(7)
    $$
    where $f = f(x,\xi,\eta),\ \lambda = \lambda(x), \  \sigma = \sigma(\xi,\eta), \
    \tau = \tau(\xi,\eta)$, neither $\lambda$ nor $\sigma^2 + \tau^2$ being equal to zero.

    In the system (7), we shall perform the following change of the variables:
    $$
    y = \omega(x), \ \mu = \varphi(\xi,\eta), \ \nu = \psi(\xi,\eta),
    $$
    which yields the following writing of its former equation:
    $$
    \lambda_1y_xf_y + (\sigma_1\mu_\xi + \tau_1\mu_\eta)f\mu +
    (\sigma_1\nu_\xi + \tau_1\nu_\eta)f_\nu = 0,
    $$
    and borrow the functions of that change from the solutions of the integrable equations
    $$
    \lambda_1\omega_x = 1, \  \sigma_1\varphi_\xi + \tau_1\varphi_\eta = 1, \
    \sigma_1\psi_\xi + \tau_1\psi_\eta = 0.
    $$

    If in the system (7) we get back to the previous designation of the coefficients and
    of the new variables, i.e. of the new coordinates in the manifolds $\mathfrak{M}$ and $\mathfrak{N}$,
    then its equations may be written in a simpler form:
    $$
    f_x + f_\xi = 0, \ \ \lambda f_x +\sigma f_\xi + \tau f_\eta = 0,
    \eqno(7')
    $$
    the inferior index "2" in the latter equation dropped as redundant.

    Let us demonstrate that in the system $(7')$ $\tau \neq 0$ and $\lambda \neq const$.
    If $\tau =0$ then from the system $(7')$ whose rank is equal to 2 we shall get $f_x = 0, \ f_\xi = 0,$ which is obviously in contradiction with
    the conditions (2). And in case $\lambda = a = const$ there appears from the system $(7')$ the connection of the derivatives $(\sigma - a)f_\xi + \tau f_\eta = 0,$ which is also in contradiction with the conditions (2).

    Let us differentiate the second equation of the system $(7')$ with respect to the variables $x$ and $\xi$,
    and then sum the results up, with consideration of the first equation;
    $$
    \lambda_xf_x + \sigma_\xi f_\xi + \tau_\xi f_\eta = 0.
    \eqno(8)
    $$

    If the rank of the augmented system $(7') + (8)$ of three equations is equal to three, then we arrive at the metric function
               being degenerated in all the coordinates,
    as $f_x = 0, \ f_\xi = 0, \ f_\eta = 0.$ So,
    $$
    \left|
    \begin{array}{ccc}
    1 & 1 & 0 \\
    \lambda & \sigma & \tau \\
    \lambda_x & \sigma_\xi & \tau_\xi
    \end{array}
    \right| =0,
    $$
    hence, by expanding the determinant, we establish that
    $$
    \lambda_x = \lambda \frac{\tau_\xi}{\tau} -
    \frac{\sigma \tau_\xi - \tau\sigma_\xi}{\tau}.
    \eqno(9)
    $$

    We shall differentiate the relation (9) with respect to the variable $x$, separating the variables:
    $$
    \frac{\lambda_{xx}}{\lambda_x} = \frac{\tau_\xi}{\tau} = a = const,
    $$
    wherefrom, after partial integration of the result we get the equations
    for the coefficients of the system  $(7'):$
    $$
    \lambda_x = a\lambda + b, \ \sigma_\xi = a\sigma + b, \ \tau_\xi = a\tau,
    \eqno(10)
    $$
    $a^2 + b^2$ being unequal to $0$, since $\lambda \neq const.$

    When integrating the equations (10) we shall consider two cases separately:
    $a = 0$ and $a \neq 0$. In the former case, that of
    $$
    a = 0.
    \eqno(11)
    $$

    the equations (10) become simpler:
    $$
    \lambda_x = b, \ \sigma_\xi = b, \ \tau_\xi = 0
    $$
    and their solutions may be written as follows:
    $$
    \lambda(x) = bx +c, \ \sigma(\xi,\eta) = b\xi + b\sigma(\eta), \
    \tau(\xi,\eta) = b\tau(\eta),
    $$
    where, obviously, $b \neq 0$ and $\tau(\eta) \neq 0.$

    Let us write down the system $(7')$ for the metric function with the coefficients found in case (11):
    $$
    f_x + f_\xi = 0, \ xf_x + (\xi + \sigma(\eta))f_\xi +
    \tau(\eta)f_\eta = 0.
    \eqno(12)
    $$

    We shall perform in the system (12) the following change of the variables:
    $$
    \mu = \xi + \varphi(\eta), \ \nu = \psi(\eta),
    \eqno(13)
    $$
    where $\psi' \neq 0,$ taking into account that with such change
    $$
    f_\xi = f_\mu, \ f_\eta = \varphi'(\eta)f_\mu + \psi'(\eta)f_\nu.
    $$

    The former of the equations (12) will retain its simplest form,
    while the latter changes:
    $$
    f_x + f_\mu = 0, \ xf_x + (\xi +\sigma(\eta) +
    \tau(\eta)\varphi'(\eta))f_\mu + \tau(\eta)\psi'(\eta)f_\nu = 0.
    $$

    The functions $\varphi$ and $\psi$ in the change of the variables (13) may be taken
    from the solutions of the equations
    $$
    \sigma(\eta) + \tau(\eta)\varphi'(\eta) = \varphi(\eta), \
    \tau(\eta)\psi'(\eta) = 1.
    $$
    If we now get back to the previous designation of the variables the system (12)is written in the following form, simplest for the case (11):
    $$
    f_x + f_\xi = 0, \ xf_x + \xi f_\xi + f_\eta = 0.
    \eqno(12')
    $$

    We shall substitute the solution $f = \theta(x + \xi,\eta)$ of the first equation of the system $(12')$, where
    $\theta = \theta(u,v)$ is a function of two variables
    $u = x + \xi$ and $v = \eta$, into its second equation: $u\theta_u +
    \theta_v = 0.$ By integrating the resulting equation, we find:
    $\theta(u,v) = \chi(u\exp{(-v)}),$ where $\chi$ is an arbitrary function of one
    variable with one derivative unequal to zero. As result,
    we get for the metric function the expression $f = \chi((x - \xi)\exp{(-\eta)}),$ which we shall rewrite
    as follows: $\chi^{-1}(f) = x\exp{(-\eta)} - \xi\exp{(-\eta)},$
    where $\chi^{-1}$ is the function inverse to $\chi$. It can be seen that with
    an accuracy up to the change of coordinates $\exp{(-\eta)} \to \xi, \
    -\xi\exp{(-\eta)} \to \eta$  and a scaling transformation $\chi^{-1}(f) \to f$
    it coincides with the canonical expression (5) that took place in the theorem
    that we proved above. We shall note that all the above said changes of variables
    performed while solving the initial system (7) are essentially changes of coordinates in the manifolds $\mathfrak{M}$ and $\mathfrak{N}$.

    Let us now consider the latter case, i.e. that when
    $$
    a \neq 0.
    \eqno(14)
    $$

    The solutions of the differential equations (10) for the coefficients of the system $(7')$ may be written down as follows:
    $$
    \lambda(x) = c\exp{ax} - b/a,\ \sigma(\xi,\eta) =
    c\sigma(\eta)\exp{a\xi} - b/a, \ \tau(\xi,\eta) = ac\tau(\eta)\exp{a\xi},
    $$
    where $c \neq 0, \ \tau(\eta) \neq 0,$  as $\lambda' \neq 0$ and
    $\tau(\xi,\eta) \neq 0$.

    We shall introduce new coefficients into the system $(7')$, introducing additionally change of the variables $ ax \to x, \ a\xi \to \xi:$
    $$
    f_x + f_\xi = 0, \ \exp{(x-\xi)}f_x + \sigma(\eta)f_\xi + \tau(\eta)f_\eta =
    0.
    \eqno(15)
    $$

    We shall substitute the solution $ f = \theta(x - \xi,\eta)$ of the first equation of the system (15) into its second equation: $(\exp u - \sigma(v))\theta_u + \tau(v)\theta_v = 0,$
    where, as in the former case, $u = x - \xi, \ v = \eta.$ The equation is solved by
    the method of characteristics:
    $$
    \theta(u,v) = \chi(\varphi(v)\exp{(-u)} + \psi(v)),
    $$
    Where the functions $\varphi(v)$ and $\psi(v)$ are, in their turn, solutions of other equations:
    $$
    \frac{\varphi'(v)}{\varphi(v)} = - \frac{\sigma(v)}{\tau(v)}, \ \
    \frac{\psi'(v)}{\varphi(v)} = \frac{1}{\tau(v)}.
    $$

    Then, for the metric function (1) we have:
         $$
    f = \chi(\varphi(\eta)\exp{(-x + \xi)} + \psi(\eta)),
    $$
    which, obviously, with an accuracy up to its scaling transformation $\chi^{-1}(f) \to f$ and change of the variables $\exp{(-x)} \to x$ and
    $\varphi(\eta)\exp \xi \to \xi, \ \psi(\eta) \to \eta$ (change of
    the coordinates in the manifolds
    $\mathfrak{M}$ and $\mathfrak{N}$) may be written in the canonical form (5).

    Let us make sure that for the metric function (5), satisfying obviously the
    inequalities (2) of the axiom {\bf A.3}, the rank of the functional
    matrix (4) is equal to 5, i.e. that for that function there exists an only
    independent equation (3) that expresses the phenomenological symmetry of rank (3,2)
    of the geometry of two sets that it defines on a one-dimensional and a two-dimensional manifolds.

    Let us substitute the metric function (5) into the functional matrix (4):
    $$
    \left\|
    \begin{array}{cccccc}
    \xi_\alpha  & \xi_\beta  & 0 & 0 & 0 & 0 \\
    0 & 0 & \xi_\alpha  & \xi_\beta  & 0 & 0 \\
    0 & 0 & 0 & 0 & \xi_\alpha  & \xi_\beta  \\
    x_i & 0 & x_j & 0 & x_k & 0 \\
    1 & 0 & 1 & 0 & 1 & 0 \\
    0 & x_i & 0 & x_j & 0 & x_k \\
    0 & 1 & 0 & 1 & 0 & 1
    \end{array}
    \right\|
    \eqno(16)
    $$
    Since the Jacobian of the fifth order of the matrix (16), obtained by
    crossing out the last two rows and the last column and equal to
    $- \xi_\alpha\xi_\beta^2(x_i - x_j)$, is unequal to zero, and any Jacobian of the sixth order is equal to zero, the rank of the matrix is equal to 5.

    For the six functions
    $$
    \left.\begin{array}{c}
    f(i\alpha)=x_i\xi_\alpha + \eta_\alpha, \
    f(i\beta)=x_i\xi_\beta + \eta_\beta, \\
    f(j\alpha)=x_j\xi_\alpha + \eta_\alpha, \
    f(j\beta)=x_j\xi_\beta + \eta_\beta, \\
    f(k\alpha)=x_k\xi_\alpha + \eta_\alpha, \
    f(k\beta)=x_k\xi_\beta + \eta_\beta
    \end{array}\right\}
    $$
    the equation (3) is found by way of elimination of all the seven coordinates
    $x_i, x_j, x_k, \xi_\alpha, \eta_\alpha, \xi_\beta, \eta_\beta$
    of the points of the cortege $<ijk,\alpha\beta>$:
    $$
    \left|
    \begin{array}{ccc}
    f(i\alpha ) & f(i\beta ) & 1 \\
    f(j\alpha ) & f(j\beta ) & 1 \\

    f(k\alpha ) & f(k\beta ) & 1
    \end{array}
    \right| =0,
    $$

    The proof of the theorem is complete.

    \vspace{5mm}

    We shall note that the PS G2S of rank (3,2) is endowed with a group symmetry of degree 2. \
    Its set of motions consists of two-parameter groups of transformations
    $$
    x' = ax + b, \ \ \xi' = \xi/a, \ \eta' = \eta - b\xi/a,
    $$
    preserving the metric function $x'\xi' + \eta' = x\xi + \eta$, i.e. the distance between any two of its points.
    Thus, the Erlangen programme of F. Klein (1872) that was formulated by him for ordinary geometries on one set is also valid for the phenomenologically symmetric geometry of two sets of rank (3,2) in question.

    We shall also note that the PS G2S of rank (3,2) springs up in a natural way when we analyze the structure of Ohm's law [1. pp. 121-123] and so can be called a physical structure.
    For the first time its rigorous mathematical definition was given in note [2],
    but the proof of the corresponding theorem of existence and singleness was more complicated and bulkier.

    \vspace{10mm}

    \begin{center}
    {\bf Литература}
    \end{center}

    \vspace{5mm}

    1. Кулаков Ю.И. Теория физических структур. - М.: Доминико, 2004. - 847 с.

    2. G.G. Mihailichenko. {\bf A binary physical structure of rank (3,2).}
    Sibarian Mathematical Journal, 1973, Volume 14, Issue 5, pp 737-742
    (The Russion version: Михайличенко Г.Г. Бинарная физическая структура ранга (3,2) // Сиб. мат.
    журн., 1973, Т.14, №5, С. 1057-1064).

    \vspace{50mm}

    { \bf Г.Г. Михайличенко} \\ профессор кафедры физики.
    \\ Горно-Алтайский государственный
    университет, \\ 649000, г. Горно-Алтайск, ул. Ленкина, д. 1, \\
    e-mail: mikhailichenko@gasu.ru

    \vspace{5mm}

    {\bf G.G. Mihailichenko} \\
    Professor of Chair of Physics, \\
    Gorno-Altaisk State University, \\
    1 Lenkin str., Gorno-Altaisk, 649000, Russia, \\
    e-mail: mikhailichenko@gasu.ru

    \vspace{200mm}

    \begin{flushright}
    УДК 514.1
    \end{flushright}

    \begin{center}
    {\Large \bf The Phenomenologically symmetric geometry \\ of two sets of
    rank (3,2)}
    \end{center}

    \begin{center}
    Mihailichenko G.G.
    \end{center}

    A geometry of two sets (GTS) is given on manifolds $\mathfrak{M}$ and
    $\mathfrak{N}$ by a metric (two-point) function
    $f:\mathfrak{M\times N}\to R$. Its
    phenomenological symmetry (PS) means that for some numbers of points
    from each manifold all the reciprocal distances are tied to some equation.
    Such simplest geometry on one-dimensional manifolds
    was discovered by Yu.I. Kulakov when he was analyzing
    the structure of Newton's 2nd law.
    In this note the PS G2S of rank (3,2) that springs up in the process of analysis the structure of Ohm's law is precisely defined and analyzed using a new approach.

    \vspace{5mm}

    Bibliography: 2 names.

    \vspace{5mm}
    [Gorno-Altaisk State University (GASU)]

    \vspace {5mm}

    { \bf Gennadiy Grigоrievich Mihailichenko,} 1942, GASU, Professor of
    Chair of Physics, Doctor of Physical and Mathematical Sciences, Professor,
    649000, Gorno-Altaisk, 1 Lenkina Street, GASU
    (649000, Gorno-Altaisk, 16 Ulagasheva Street, flat 11), {\bf kfizika@gasu.ru} \\
    {\bf (mikhailichenko@gasu.ru)}. 8-388-22-2-75-39 (8-388-22-2-86-31).

    \end{document}